\documentclass{amsart}

\usepackage{amssymb,latexsym,amsfonts,amsmath}
\usepackage{graphicx}
\include{diagrams}
\usepackage{eso-pic}
\usepackage{graphicx}
\usepackage{color}

\usepackage{amssymb,latexsym,amsfonts,amsmath}
\usepackage{graphicx}

\newtheorem{theorem}{Theorem}[section]
\newtheorem{lemma}[theorem]{Lemma}

\newtheorem{proposition}[theorem]{Proposition}
\newtheorem{corollary}[theorem]{Corollary}

\theoremstyle{definition}
\newtheorem{definition}[theorem]{Definition}
\newtheorem{example}[theorem]{Example}

\theoremstyle{remark}

\numberwithin{equation}{section}

\renewcommand{\Re}{{\mathbb{R}}}
\newcommand{\R}{{\mathbb{R}}}

\newcommand{\kr}{\textrm{ker}}

\newcommand{\K}{\mathcal{K}}
\newcommand{\Kinfty}{\K_{\infty}}
\newcommand{\x}{\mathbf{x}}
\newcommand{\y}{\mathbf{y}}
\newcommand{\xu}{\mathbf{x_u}}
\renewcommand{\u}{\mathbf{u}}

\begin{document}
\begin{abstract}
The reduction of dynamical systems has a rich history, with many
important applications related to stability, control and
verification.  Reduction of nonlinear systems is typically performed in an ``exact''
manner---as is the case with mechanical systems with
symmetry---which, unfortunately, limits the type of systems to which
it can be applied.  The goal of this paper is to consider a more
general form of reduction, termed \emph{approximate reduction}, in
order to extend the class of systems that can be reduced.  Using
notions related to incremental stability, we give conditions on when
a dynamical system can be projected to a lower dimensional space
while providing hard bounds on the induced errors, i.e., when it is
behaviorally similar to a dynamical system on a lower dimensional
space.  These concepts are illustrated on a series of examples.
\end{abstract}

\title[Approximate Reduction of Dynamical Systems]{Approximate Reduction of Dynamical Systems}
\thanks{This research was partially supported by the National Science Foundation, EHS award 0712502.}

\author[Paulo Tabuada]{Paulo Tabuada}
\address{Department of Electrical Engineering\\
University of California at Los Angeles,\\
Los Angeles, CA 90095}
\email{tabuada@ee.ucla.edu}

\author[Aaron D. Ames]{Aaron D. Ames}
\address{Control and Dynamical Systems Department\\
California Institute of Technology\\
Pasadena, CA 91125}
\email{ames@cds.caltech.edu}

\author[Agung Julius]{Agung Julius}
\address{Department of Electrical and Systems Engineering\\
University of Pennsylvania\\
Philadelphia, PA 19104}
\email{agung@seas.upenn.edu}

\author[George J. Pappas]{George J. Pappas}
\address{Department of Electrical and Systems Engineering\\
University of Pennsylvania\\
Philadelphia, PA 19104}
\email{pappasg@ee.upenn.edu}

\maketitle

\section{Introduction}


Modeling is an essential part of many engineering disciplines and
often a key ingredient for successful designs. Although it is
widely recognized that models are only approximate descriptions of
reality, their value lies precisely on the ability to describe,
within certain bounds, the modeled phenomena. In this paper we
consider modeling of closed-loop nonlinear control systems, i.e.,
differential equations, with the purpose of simplifying the
analysis of these systems. The goal of this paper is to reduce the
dimensionality of the differential equations being analyzed while
providing hard bounds on the introduced errors.  One promising
application of these techniques is to the verification of hybrid
systems, which is currently constrained by the complexity of high
dimensional differential equations.

Reducing differential equations---and in particular mechanical
systems---is a subject with a long and rich history.  The first form
of reduction was discovered by Routh in the 1860's; over the years,
geometric reduction has become an academic field in itself.  One
begins with a differential equation with certain symmetries, i.e.,
it is invariant under the action of a Lie group on the phase space.
Using these symmetries, one can reduce the dimensionality of the
phase space (by ``dividing'' out by the symmetry group) and define a
corresponding differential equation on this reduced phase space. The
main result of geometric reduction is that one can understand the
behavior of the full-order system in terms of the behavior of the
reduced system and vice versa~\cite{SymplecticReduction,Noether,NonH}.  While this form of ``exact'' reduction is very
elegant, the class of systems for which this procedure can be
applied is actually quite small.  This indicates the need for a form
of reduction that is applicable to a wider class of systems and,
while not being exact, is ``close enough''.

In systems theory, reduced order modeling has also been
extensively studied under the name of model reduction~\cite{MReduc,SurveyModelReduction}. Contrary to
model reduction where approximation is measured using $L^2$ norms
we are interested in $L^\infty$ norms. The guarantees provided by $L^\infty$ norms are more natural when applications to safety verification are of interest. More recent work considered the exact reduction of control systems~\cite{BisimTAC,BisimSCL} based on the notion of bisimulation which was later generalized to approximate bisimulation~\cite{ApproxBisimDyn,Tab07,ApproxBisimLinear}.

We develop our results in the framework of incremental stability
and our main result is in the spirit of existing stability results
for cascade systems that proliferate the Input-to-State Stability
(ISS) literature. See for example~\cite{ISS} and the references therein. A preliminary version of our results appeared in the conference paper~\cite{TAJPCDC07}.


\section{Preliminaries}

A continuous function $\gamma:\Re_0^+\to \Re_0^+$, is said to belong
to class $\Kinfty$ if it is strictly increasing, $\gamma(0)=0$ and
$\gamma(r)\to \infty$ as $r\to\infty$. A continuous function
$\beta:\Re_0^+\times\Re_0^+\to\Re_0^+$ is said to belong to class
$\mathcal{KL}$ if, for each fixed $s$, the map $\beta(r,s)$ belongs
to class $\mathcal{K}_\infty$ with respect to $r$ and, for each
fixed $r$, the map $\beta(r,s)$ is decreasing with respect to $s$
and $\beta(r,s)\to 0$ as $s\to\infty$.

For a smooth function $\phi:\Re^n\to\Re^m$ we denote by $T\phi$ the
tangent map to $\phi$ and by $T_x\phi$ the tangent map to $\phi$ at
$x\in \Re^n$. We will say that $\phi$ is a submersion at $x\in
\Re^n$ if $T_x\phi$ is surjective and that $T\phi$ is a submersion
if it is a submersion at every $x\in\Re^n$. When $\phi$ is a
submersion we will also use the notation $\kr(T\phi)$ to denote the
distribution:
$$\kr(T\phi)=\{X:\Re^n\to\Re^n\,\,\vert\,\, T\phi\cdot X=0\}.$$
The Lie bracket of vector fields $X$ and $Y$ will be denotes by $[X,Y]$.

Given a point $x\in \Re^n$, $\vert x\vert$ will denote the usual
Euclidean norm while $\vert\vert f\vert\vert$ will denote
$\mathrm{ess}\sup_{t\in [0,\tau]}\vert f(t)\vert$ for any given
function $f:[0,\tau]\to \Re^n$, $\tau\in\R^+$.

\subsection{Dynamical and control systems}
In this paper we shall restrict our attention to dynamical and
control systems defined on Euclidean spaces.

\begin{definition}
A {\it vector field} is a pair $(\Re^n,X)$ where $X$ is a smooth map
\mbox{$X:\Re^n\to \Re^n$}. A smooth curve
$\mathbf{x}(\,\cdot\,,x):I\to \Re^n$, defined on an open subset $I$
of $\R$ including the origin, is said to be a trajectory of
$(\Re^n,X)$ if the following two conditions hold:
\begin{enumerate}
\item $\mathbf{x}(0,x)=x$;
\item $\frac{d}{dt}\mathbf{x}(t,x)=X(\x(t,x))$ for all $t\in I$.
\end{enumerate}
\end{definition}

A control system can be seen as an under-determined vector field.

\begin{definition}
A {\it control system} is a triple $(\Re^n,\Re^m,F)$ where $F$ is a
smooth map \mbox{$F:\Re^n\times\Re^m\to \Re^n$}. A smooth curve
$\mathbf{x_u}(\cdot,x):I\to \Re^n$, defined on an open subset $I$ of
$\R$ including the origin, is said to be a trajectory of
$(\Re^n,\Re^m,F)$ if there exists a smooth curve
$\mathbf{u}:I\to\Re^m$ such that the following two conditions hold:
\begin{enumerate}
\item $\mathbf{x_u}(0,x)=x$;
\item $\frac{d}{dt}\mathbf{x_u}(t,x)=F(\xu(t,x),\u(t))$ for all $t\in I$.
\end{enumerate}
\end{definition}

We have defined trajectories based on smooth input curves mainly for simplicity since the presented results hold under weaker regularity assumptions.

 
\section{Exact reduction}

For some dynamical systems described by a vector field $X$ on $\Re^n$
it is possible to replace $X$ by a vector field $Y$ describing the
dynamics of the system on a lower dimensional space, $\Re^m$, while
retaining much of the information in $X$. When this is the case we
say that $X$ can be reduced to $Y$. This idea of (exact) reduction
is captured by the notion of $\phi$-related vector fields.

\begin{definition}
\label{FRel} Let $\phi:\Re^n\to\Re^m$ be a smooth map. The vector
field $(\Re^n,X)$ is said to be {\it $\phi$-related} to the vector
field $(\Re^m,Y)$ if for every $x\in\R^n$ we have:
\begin{equation}
T_x\phi\cdot X(x)=Y\circ \phi(x).
\end{equation}
\end{definition}

The following proposition, proved in~\cite{MTAA}, characterizes
$\phi$-related vector fields in terms of their trajectories.

\begin{proposition}
{\it The vector field $(\Re^n,X)$ is $\phi$-related to the vector
field $(\Re^m,Y)$ for some smooth map \mbox{$\phi:\Re^n\to\Re^m$} iff for every $x\in\R^n$ we have:
\begin{equation}
\label{ExactReduc} \phi\circ \x(t,x)=\y(t,\phi(x)),
\end{equation}
where $\x(t,x)$ and $\y(t,y)$ are the trajectories of vector fields
$X$ and $Y$, respectively.
}
\end{proposition}

For $\phi$-related vector fields, we can replace the study of
trajectories $\x(\cdot,x)$ with the study of trajectories
$\y(\cdot,\phi(x))$ living on the lower dimension space $\Re^m$. In particular, formal verification of $X$ can be performed on $Y$
whenever the relevant sets describing the verification problem can
also be reduced to $\Re^m$.

If a vector field and a submersion $\phi$ are given we can use the
following result, proved in~\cite{Dynamics}, to determine the
existence of $\phi$-related vector fields.

\begin{proposition}
\label{Existence} {\it Let $(\Re^n,X)$ be a vector field and let
$\phi:\Re^n\to\Re^m$ be a smooth submersion. There exists a vector
field $(\Re^m,Y)$ that is $\phi$-related to $(\Re^n,X)$ iff:
$$[\kr(T\phi),X]\subseteq \kr(T\phi).$$}
\end{proposition}

In an attempt to enlarge the class of vector fields that can be
reduced we introduce, in the next section, an approximate notion of
reduction.

 
\section{Approximate Reduction}

The generalization of Definition~\ref{FRel} proposed in this section requires a decomposition of
$\Re^n$ of the form \mbox{$\Re^n=\Re^m\times\Re^k$.} Associated with this
decomposition are the canonical projections  $\pi_m:\Re^n\to\Re^m$
and $\pi_k:\Re^n\to\Re^k$ taking $\Re^n\ni
x=(y,z)\in\Re^m\times\Re^k$ to $\pi_m(x)=y$ and $\pi_k(x)=z$,
respectively.

\begin{definition}
\label{App} The vector field $(\Re^n,X)$ is said to be {\it
approximately $\pi_m$-related} to the vector field $(\Re^m,Y)$ if
there exist a class $\mathcal{K}_\infty$ function $\gamma$ and a constant $c\in\R_0^+$ such that the following estimate holds for every $x\in\R^n$:
\begin{equation}
\label{AIneq} \vert\pi_m\circ
\x(t,x)-\y(t,\pi_m(x))\vert\le\gamma(\vert\pi_k(x)\vert)+c.
\end{equation}
\end{definition}

Note that when $X$ and $Y$ are $\pi_m$-related we have:
$$\vert\pi_m\circ \x(t,x)-\y(t,\pi_m(x))\vert=0,$$
which implies (\ref{AIneq}). Definition~\ref{App} can thus be seen
as a generalization of exact reduction captured by Definition~\ref{FRel}. Similar ideas have been used in the context of approximate notions of equivalence for control systems~\cite{ApproxBisimLinear}.

Although the bound on the gap between the projection of the
original trajectory $\x$ and the trajectory $\y$ of the
approximate reduced system is a function of $x$, in concrete applications
the initial conditions are typically restricted to a bounded set
of interest. The following result has interesting implications in
these situations.

\begin{proposition}
If $(\Re^n,X)$ is approximately $\pi_m$-related to $(\Re^m,Y)$
then for any compact set $S\subseteq \Re^n$ there exists a
$\delta\in\R^+$ such that for all $x\in S$ the following estimate
holds:
\begin{equation}
\vert\pi_m\circ \x(t,x)-\y(t,\pi_m(x))\vert\le\delta.
\end{equation}
\end{proposition}

{\bf Proof:}
Let $\delta=\max_{x\in C}\gamma(\vert\pi_k(x)\vert)+c$. The scalar
$\delta$ is well defined since
$\gamma(\vert\pi_k(\:\cdot\:)\vert)+c$ is a continuous map and $C$
is compact.
$\qed$

From a practical point of view, approximate reduction is only a
useful concept if it admits characterizations that are simple to
check. In order to derive such characterizations we need to review
several notions of incremental stability.

 
\subsection{Incremental stability}

In this subsection we review two notions of incremental stability
which will be fundamental in proving the main contribution of this
paper. We follow~\cite{UBIBS} and~\cite{IncrementalS}.

\begin{definition}
A control system $(\Re^n,\Re^m,F)$ is said to be {\it incrementally
uniformly bounded-input-bounded-state stable (IUBIBSS)} if there
exist two class $\mathcal{K}_\infty$ functions $\gamma_1$ and
$\gamma_2$ such that for each $x_1,x_2\in \Re^n$ and for each pair
of smooth curves \mbox{$\mathbf{u}_1,\mathbf{u}_2:I\to\Re^m$} the following
estimate holds:
\begin{equation}
\label{EstIUBIBSS}
\vert\mathbf{x_{u_1}}(t,x_1)-\mathbf{x_{u_2}}(t,x_2)\vert\le\gamma_1(\vert
x_1-x_2\vert)+\gamma_2(\Vert\mathbf{u}_1-\mathbf{u}_2\Vert)
\end{equation}
for all $t\in I$.
\end{definition}

In general it is difficult to establish IUBIBSS directly. A
sufficient condition is given by the existence of an IUBIBSS
Lyapunov function. Note, however, that IUBIBSS only implies the
existence of a IUBIBSS Lyapunov function with very weak regularity
conditions~\cite{UBIBS}.

\begin{definition}
\label{IUBIBSSL} A $C^1$ function $V:\Re^n\times\Re^n\to\Re_0^+$ is
said to be an {\it IUBIBSS Lyapunov function} for control system
$(\Re^n,\Re^m,F)$ if there exist a $\xi\in\R^+$ and class
$\mathcal{K}_\infty$ functions
$\underline{\alpha},\overline{\alpha}$, and $\mu$  such that for every $x_1,x_2\in\R^n$ and $u_1,u_2\in\R^m$ the following holds:
\begin{enumerate}
\item \mbox{$\vert x_1-x_2\vert\ge \xi\quad\implies\quad\underline{\alpha}(\vert x_1-x_2\vert)\le V(x_1,x_2)\le\overline{\alpha}(\vert x_1-x_2\vert);$}
\item $\mu(r)\ge r+\xi$ for $r\in\Re_0^+$;
\item $\vert x_1-x_2\vert\ge \mu(\vert u_1-u_2\vert)\quad\implies\quad\dot{V}\le 0$.
\end{enumerate}
\end{definition}

A stronger notion than IUBIBSS is incremental input-to-state
stability.

\begin{definition}
A control system $(\Re^n,\Re^m,F)$ is said to be {\it incrementally
input-to-state stable (IISS)} if there exist a class $\mathcal{KL}$
function $\beta$ and a class $\mathcal{K}_\infty$ function $\gamma$
such that for each $x_1,x_2\in \Re^n$ and for each pair of smooth
curves \mbox{$\mathbf{u}_1,\mathbf{u}_2:I\to\Re^m$} the following estimate
holds:
\begin{equation}
\label{EstIISS}
\vert\mathbf{x_{u_1}}(t,x_1)-\mathbf{x_{u_2}}(t,x_2)\vert\le\beta(\vert
x_1-x_2\vert,t)+\gamma(\Vert\mathbf{u}_1-\mathbf{u}_2\Vert)
\end{equation}
\end{definition}

Since $\beta$ is a decreasing function of $t$ we immediately see
that~(\ref{EstIISS}) implies~(\ref{EstIUBIBSS}) with
$\gamma_1(r)=\beta(r,0)$ and $\gamma_2(r)=\gamma(r)$,
$r\in\Re_0^+$. Once again, IISS is implied by the existence of an
IISS Lyapunov function. See~\cite{IncrementalS} for a converse
result when the inputs take values in a compact set.

\begin{definition}
A $C^1$ function $V:\Re^n\times\Re^n\to\Re_0^+$ is said to be an
{\it IISS Lyapunov function} for the control system
$(\Re^n,\Re^m,F)$ if there exist class $\mathcal{K}_\infty$
functions $\underline{\alpha},\overline{\alpha},\alpha$, and $\mu$
such that:
\begin{enumerate}
\item $\underline{\alpha}(\vert x_1-x_2\vert)\le V(x_1,x_2)\le\overline{\alpha}(\vert x_1-x_2\vert)$;
\item $\vert x_1-x_2\vert\ge \mu(\vert u_1-u_2\vert)\:\:\:\implies\:\:\: \dot{V}\le -\alpha(\vert x_1-x_2\vert)$.
\end{enumerate}
\end{definition}

 
\subsection{Fiberwise stability}

In addition to incremental stability we will also need a notion of
partial practical stability.

\begin{definition}
A vector field $(\Re^n,X)$ is said to be fiberwise practically stable with
respect to $\pi_k$ if there exist a class $\mathcal{K}_\infty$
function $\gamma$ \and a constant $c\in\R_0^+$ such that the following estimate holds:
$$
\| \pi_k(\x(\cdot,x)) \| \leq \gamma(\vert \pi_k(x)\vert)+c.
$$
\end{definition}

Fiberwise practical stability can be checked with the help of the following
result.

\begin{lemma}
A vector field $(\Re^n,X)$ is fiberwise practically stable with respect to
$\pi_k$ if there exist two $\K_\infty$ functions, $\underline{\alpha}$ and
$\overline{\alpha}$, a constant $d\in\Re_0^+$, and a function $V: \Re^n \to \Re$ such that for every $x\in\R^n$ satisfying $\vert \pi_k(x)\vert \ge d$ we have:
\begin{enumerate}
\item $\underline{\alpha}(|\pi_k(x)|) \leq V(x) \leq  \overline{\alpha}(|\pi_k(x)|)$,
\item $\dot{V} \leq 0$.
\end{enumerate}
\end{lemma}

 
\subsection{Existence of approximate reductions}

In this subsection we prove the main result providing sufficient
conditions for the existence of approximate reductions.

\begin{theorem}
\label{MainTh} {\it Let $(\Re^n,X)$ be a fiberwise practically stable vector
field with respect to $\pi_k$ and let $F = T \pi_m \cdot X : \Re^m
\times \Re^k \to \Re^m$, viewed as a control system with state
space $\Re^m$, be IUBIBSS.  Then, the vector field $(\Re^m,Y)$
defined by:
$$Y(y) = T_{(y,0)} \pi_m \cdot X(y,0) = F(y,0)$$
for every $y\in\R^m$ is approximately $\pi_m$-related to $(\Re^n,X)$.}
\end{theorem}

{\bf Proof:}
By assumption, control system $(\R^m,\R^k,F=T\pi_m\circ X)$ is IUBIBSS. If we denote by $\mathbf{y}$ a trajectory of $F$ we have:
$$\vert\mathbf{y}_{\mathbf{v}_1}(t,y_{1})-\mathbf{y}_{\mathbf{v}_2}(t,y_{2})\vert \le \gamma_1(\vert y_{1}-y_{2}\vert)+\gamma_2(\Vert \mathbf{v}_1-\mathbf{v}_2\Vert).$$
In particular, we can take:
$$y_{1}=y_{2}=\pi_m(x), \quad
\mathbf{v}_1=\pi_k\circ \x(\cdot,x), \quad \mathbf{v}_2=0,
$$
to get:
\begin{eqnarray}
\vert\pi_m\circ\mathbf{x}(t,x)-\mathbf{y}(t,\pi_m(x))\vert & = & \vert\mathbf{y}_{\pi_k\circ\x(t,x)}(t,\pi_m(x))-\mathbf{y}_{0}(t,\pi_m(x))\vert\notag\\
& = & \vert\mathbf{y}_{\mathbf{v}_1}(t,\pi_m(x))-\mathbf{y}_{0}(t,\pi_m(x))\vert\notag\\
& \le &\gamma_2(\Vert \mathbf{v}_1\Vert)=\gamma_2(\Vert\pi_k\circ
\x(\cdot,x)\Vert).\notag
\end{eqnarray}
But it follows from fiber practical stability of $X$ with respect to $\pi_k$
that:
$$
\| \pi_k\circ \x(\cdot,x) \| \le\gamma(\vert \pi_k(x)\vert)+c.
$$
We thus have:
\begin{eqnarray}
\vert\pi_m\circ\mathbf{x}(t,x)-\mathbf{y}(t,\pi_m(x))\vert & \le &
\gamma_2\big( \gamma(\vert \pi_k(x)\vert) +c\big)\notag\\
& \le &\gamma_2\big(\lambda_1\gamma(\vert\pi_k(x)\vert)\big)+\gamma_2\big(\lambda_2 c\big),\notag
\end{eqnarray}
for some constants $\lambda_1,\lambda_2\in\R_0^+$. This concludes the proof since $\gamma_2(\lambda_1\gamma(\vert\cdot\vert))$ is a class
$\mathcal{K}_\infty$ function and $\gamma_2(\lambda_2 c)\in\R_0^+$.
$\qed$

Theorem~\ref{MainTh} shows that sufficient conditions for approximate reduction can be given in terms of ISS-like Lyapunov functions and how reduced system can be constructed. Before illustrating Theorem~\ref{MainTh} with several examples in the next section we present an important corollary.

\begin{corollary}
\label{Cor} {\it Let $(\Re^n,X)$ and $(\Re^m,Y)$ be vector fields satisfying the assumptions of Theorem~\ref{MainTh}.
Then, for any compact set $S\subseteq \Re^n$ there exists a
$\delta>0$ such that for any $x\in S$ and $y\in\pi_m(S)$ the
following estimate holds:
$$\vert \pi_m\circ \x(t,x)-\y(t,y)\vert\le\delta$$}
\end{corollary}

{\bf Proof:}
Using the same proof as for Theorem~\ref{MainTh}, except picking $y_1
= \pi_m(x)$ and $y_2 = y$, it follows that:
$$
\vert \pi_m\circ \x(t,x)-\y(t,y)\vert\le \gamma_1(|\pi_m(x) - y|)
+ \gamma_2\big(\lambda_1\gamma(\vert\pi_k(x)\vert)\big)+\gamma_2\big(\lambda_2 c\big).
$$
The bound $\delta$ is now given by:
$$
\delta = \max_{(x,y) \in S \times \pi_m(S)}\,\,\gamma_1(|\pi_m(x) -
y|) + \gamma_2\big(\lambda_1\gamma(\vert\pi_k(x)\vert)\big)+\gamma_2\big(\lambda_2 c\big)
$$
which is well defined since $S\times \pi_m(S)$ is compact.
$\qed$



\section{Examples}

In this section, we consider examples that illustrate the
usefulness of approximate reduction.

\begin{figure}
\begin{center}
\includegraphics[scale=0.2]{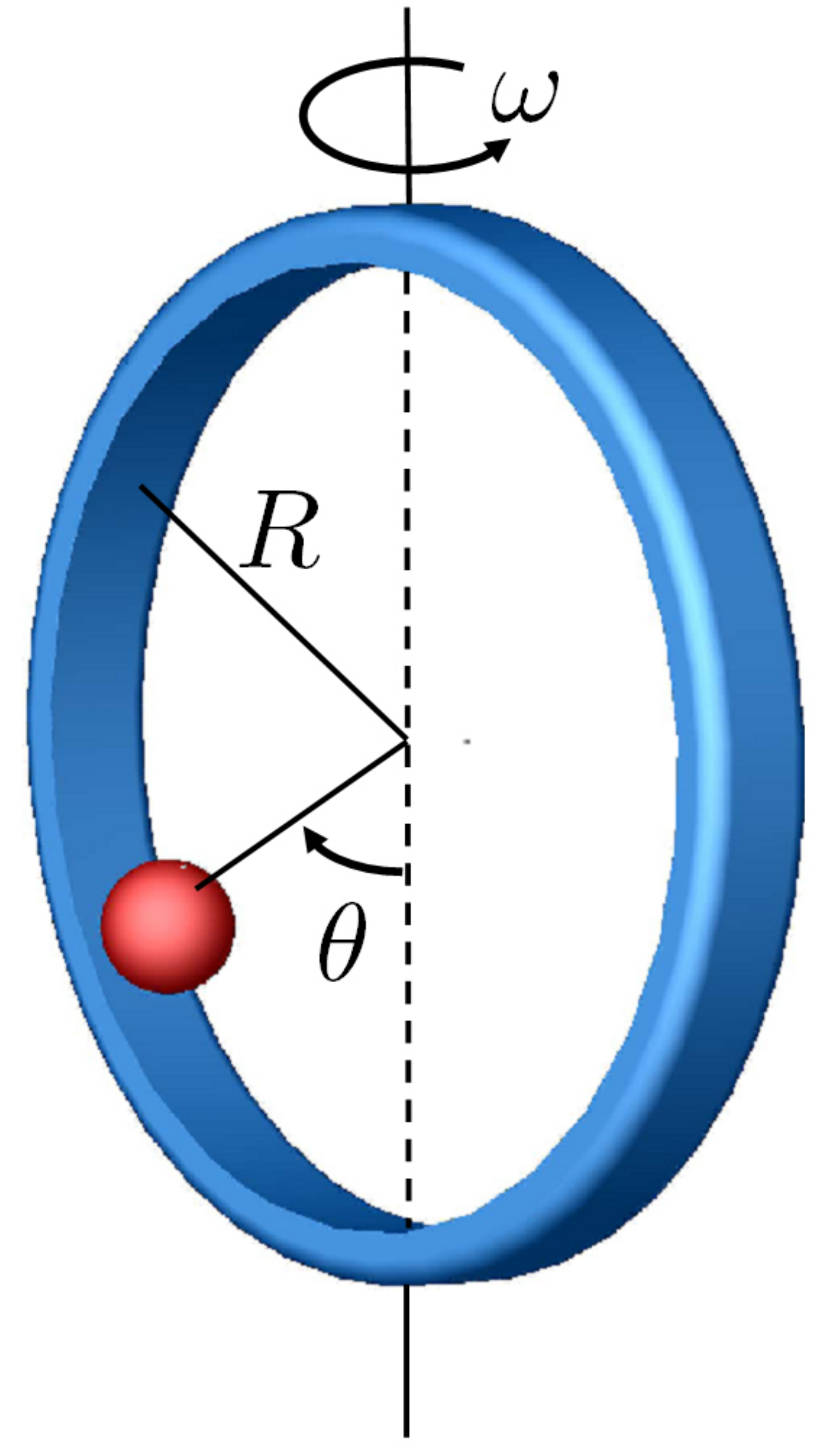}
\caption{Ball in a rotating hoop.} \label{Ball}
\end{center}
\end{figure}

\begin{figure}
\begin{center}
\includegraphics[scale=0.39]{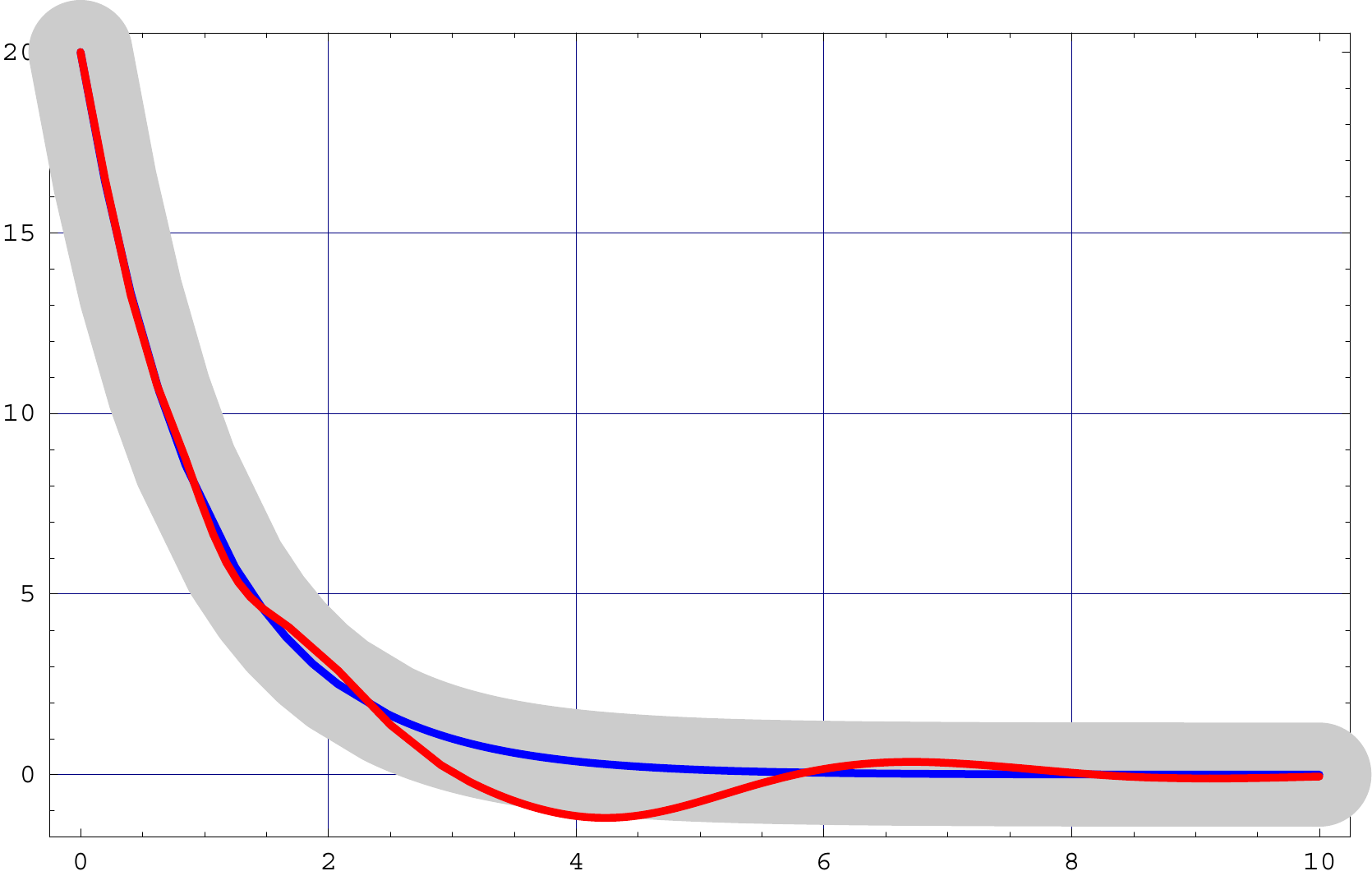}
\includegraphics[scale=0.39]{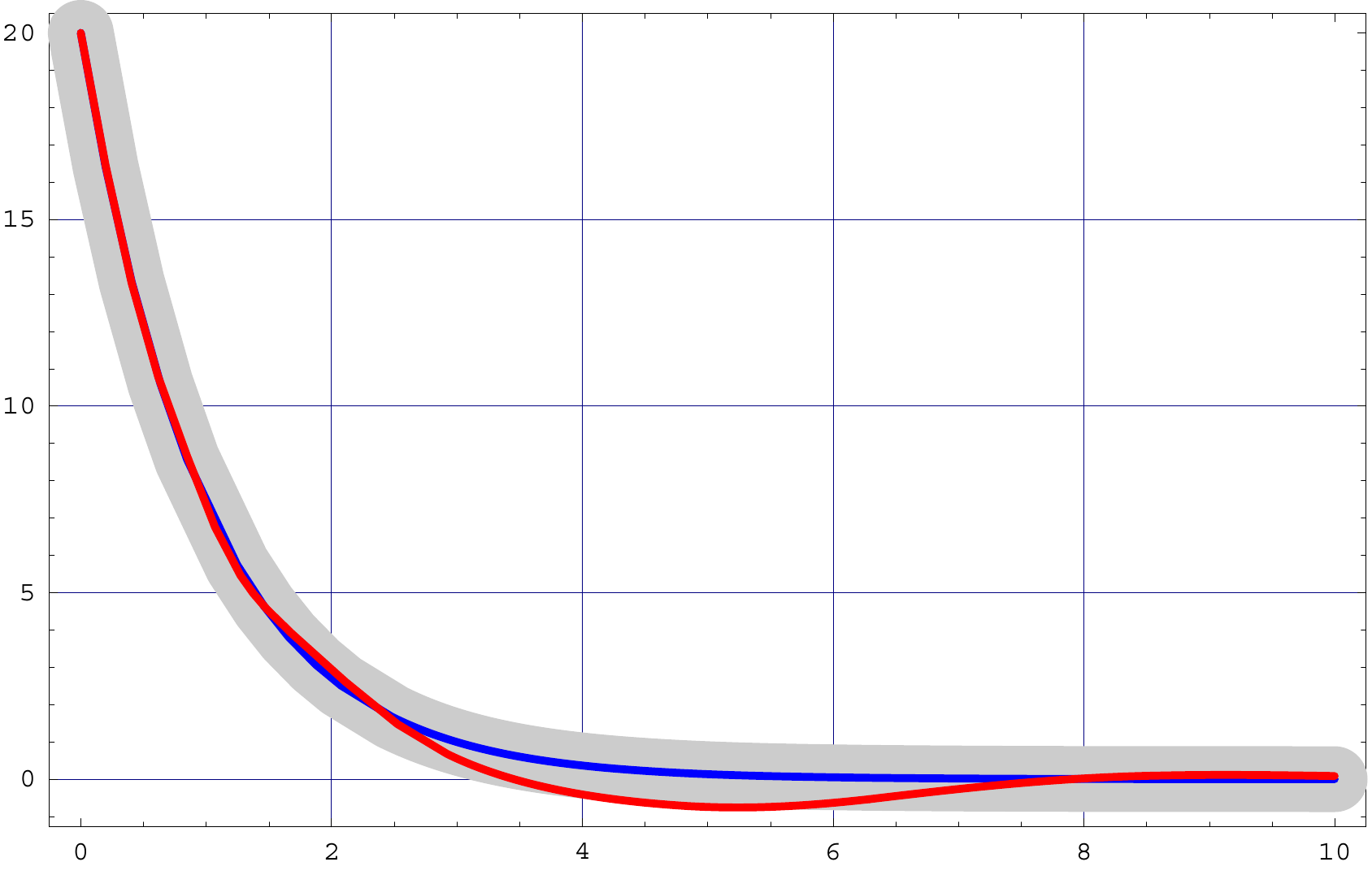}
\includegraphics[scale=0.39]{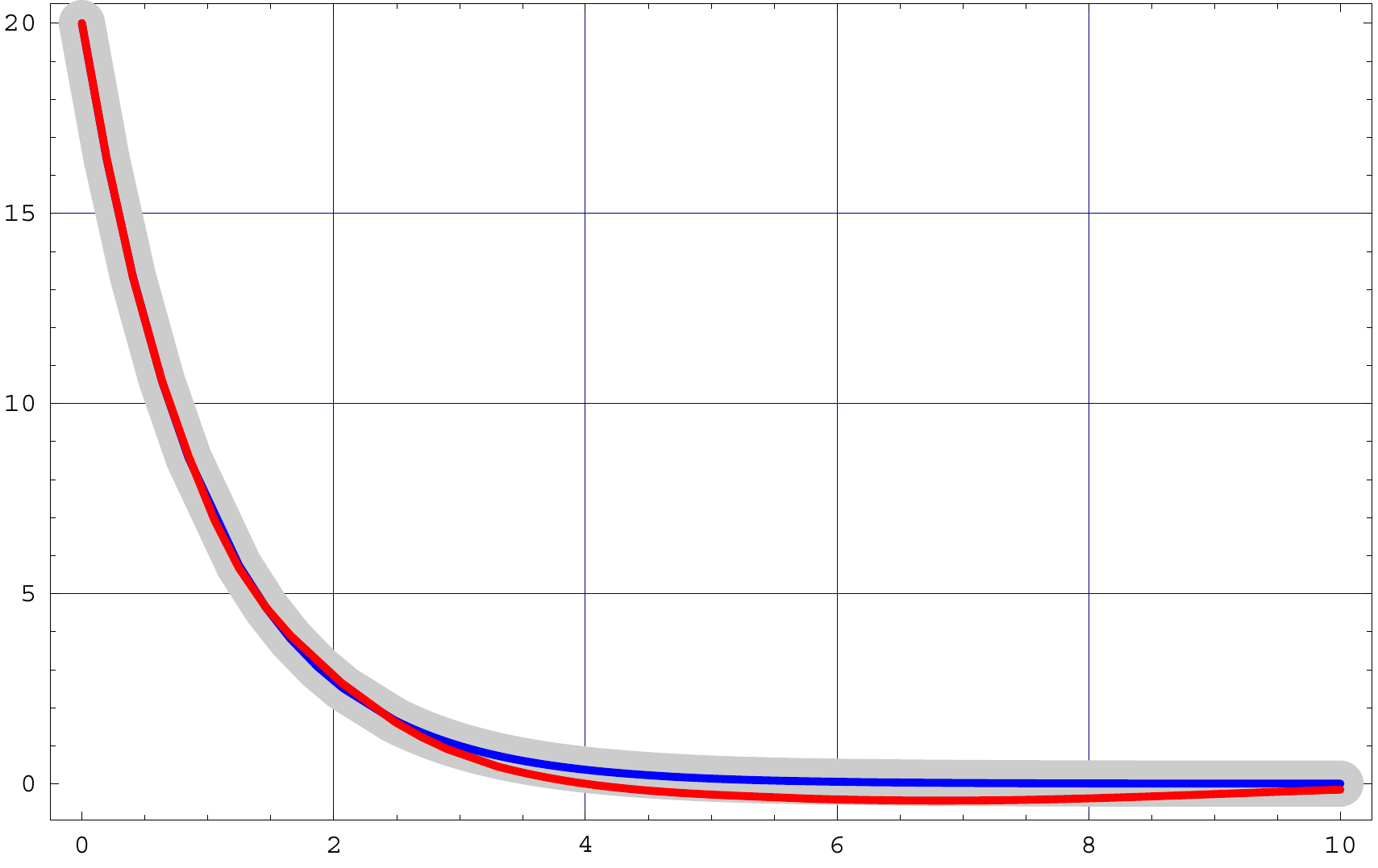}
\includegraphics[scale=0.39]{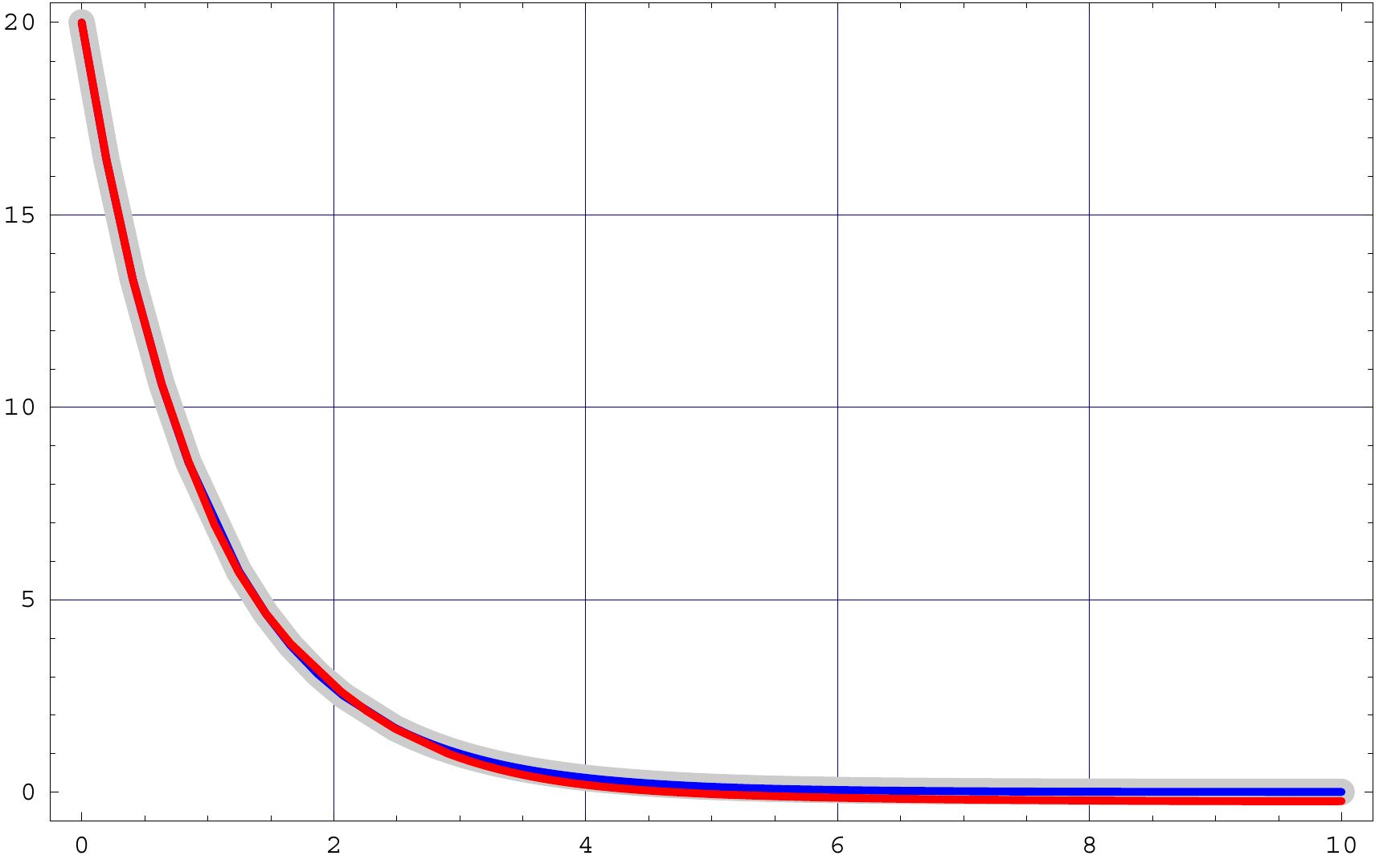}
\caption{A trajectory of the full order system (red) vs. a
trajectory for the reduced system (blue) for $R = 5,10,20,40$
(from left to right and top to bottom, respectively).} \label{example1fig}
\end{center}
\end{figure}



\begin{example}
As a first example we consider the ball in a rotating hoop with
friction, as described in Chapter 2 of~\cite{MechSym} and displayed in Figure~\ref{Ball}. For this
example, there are the following parameters:
\begin{eqnarray}
m & = & \mathrm{mass \:\: of \:\: the \:\: ball}, \nonumber\\
R & = & \mathrm{Radius \:\: of \:\: the \:\: hoop}, \nonumber\\
g & = & \mathrm{acceleration \:\: due \:\: to \:\: gravity}, \nonumber\\
\mu & = & \mathrm{friction \:\: constant \:\:for \:\: the \:\:
ball}. \nonumber
\end{eqnarray}
The equations of motion are given by:
\begin{eqnarray}
\label{BiH}
\dot{\omega} & = & -\frac{\mu}{m}\omega+\xi^2\sin\theta\cos\theta-\frac{g}{R}\sin\theta\notag\\
\dot{\theta} & = & \omega
\end{eqnarray}
where $\theta$ is the angular position of the ball and $\omega$ is
its angular velocity.

If $\pi_\omega:\Re^2\to\Re$ is the projection
$\pi_\omega(\omega,\theta)=\omega$, then according to
Proposition~\ref{Existence} there exists no vector field $Y$ on
$\Re$ which is $\pi_\omega$-related to $X$ (as defined
by~(\ref{BiH})).  However, we will show that
$Y(\omega)=T_{(\omega,0)}\pi_\omega\cdot X(\omega,0)$ is
approximate $\pi_\omega$-related to $X$.

First, we use:
$$V=\frac{1}{2}mR^2\omega^2+mgR(1-\cos\theta)-\frac{1}{2}mR^2\xi^2\sin^2\theta$$
as a Lyapunov function to show that~(\ref{BiH}) is stable. Note
that $V(\omega,\theta)=0$ for $(\omega,\theta)=(0,0)$ and
$V(\omega,\theta)>0$ for $(\omega,\theta)\ne(0,0)$ provided that
$R\xi^2<g$, which we assume. Computing the time derivative of $V$ we obtain:
$$\dot{V}=-\mu R^2\omega^2\le 0,$$
thus showing stability of~(\ref{BiH}). We now consider a compact set
$C$ invariant under the dynamics and restrict our analysis to
initial conditions in this set. Such a set can be constructed, for
example, by taking $\{x\in \Re^2\,\,\vert\,\, V(x)\le c\}$ for some
positive constant $c$. Note that stability of~(\ref{BiH}) implies
fiberwise stability on $C$ since $\pi_m(C)$ is compact.

To apply Theorem~\ref{MainTh} we only need to show that:
$$T_{(\omega,\theta)}\pi_\omega\cdot X(\omega,\theta)=-\frac{\mu}{m}\omega+\xi^2\sin\theta\cos\theta-\frac{g}{R}\sin\theta$$
is IUBIBSS on $C$ with $\theta$ seen as an input. We will conclude
IUBIBSS by proving the stronger property of IISS.
Consider the function:
$$U=\frac{1}{2}(\omega_1-\omega_2)^2.$$
Its time derivative is given by:
\begin{eqnarray}
\label{DerivU}
\dot{U} & = & (\omega_1-\omega_2)\Big[-\frac{\mu}{m}(\omega_1-\omega_2)+\xi^2\sin\theta_1\cos\theta_1\notag\\
&&-\frac{g}{R}\sin\theta_1 -\xi^2\sin\theta_2\cos\theta_2+\frac{g}{R}\sin\theta_2\Big]\notag\\
& \le & -\frac{\mu}{m}(\omega_1-\omega_2)^2+\vert \omega_1-\omega_2\vert \Big\vert\xi^2\sin\theta_1\cos\theta_1\notag\\
&& -\frac{g}{R}\sin\theta_1-\xi^2\sin\theta_2\cos\theta_2+\frac{g}{R}\sin\theta_2 \Big\vert\notag\\
& \le & -\frac{\mu}{m}(\omega_1-\omega_2)^2+\vert \omega_1-\omega_2\vert L\vert\theta_1-\theta_2\vert\notag\\
& = &
-\frac{\mu}{2m}(\omega_1-\omega_2)^2\\
&&+\Big(-\frac{\mu}{2m}(\omega_1-\omega_2)^2+\vert
\omega_1-\omega_2\vert L\vert\theta_1-\theta_2\vert\Big), \notag
\end{eqnarray}
where the second inequality follows from the fact that
$\xi^2\sin\theta\cos\theta-\frac{g}{R}\sin\theta$ is a smooth
function defined on the compact set $\pi_\theta(C)$ and is
thus globally Lipschitz on $\pi_\theta(C)$ (since its derivative is
continuous and thus bounded on any compact convex set containing $\pi_\theta(C)$) with Lipschitz
constant $L$. We now note that the condition:
$$\vert \omega_1-\omega_2\vert>\frac{2mL}{\mu}\vert \theta_1-\theta_2\vert$$
makes the second term in~(\ref{DerivU}) negative from which we
conclude the following implication:
$$\vert \omega_1-\omega_2\vert>\frac{2mL}{\mu}\vert \theta_1-\theta_2\vert\quad\implies\quad\dot{U}\le-\frac{\mu}{2m}(\omega_1-\omega_2)^2$$
showing that $U$ is an IISS Lyapunov function. We can thus
reduce~(\ref{BiH}) to:
$$\dot{\omega}=-\frac{\mu}{m}\omega.$$
Projected trajectories of the full-order system as compared with
trajectories of the reduced system can be seen in Figure
\ref{example1fig}; here $\mu = m = 1$ and $\xi = 0.1$.  Note that
as $R \to \infty$, the reduced system converges to the full-order
system (or the full-order system effectively becomes decoupled).
\end{example}

\begin{figure}
\begin{center}
\includegraphics[scale=0.3]{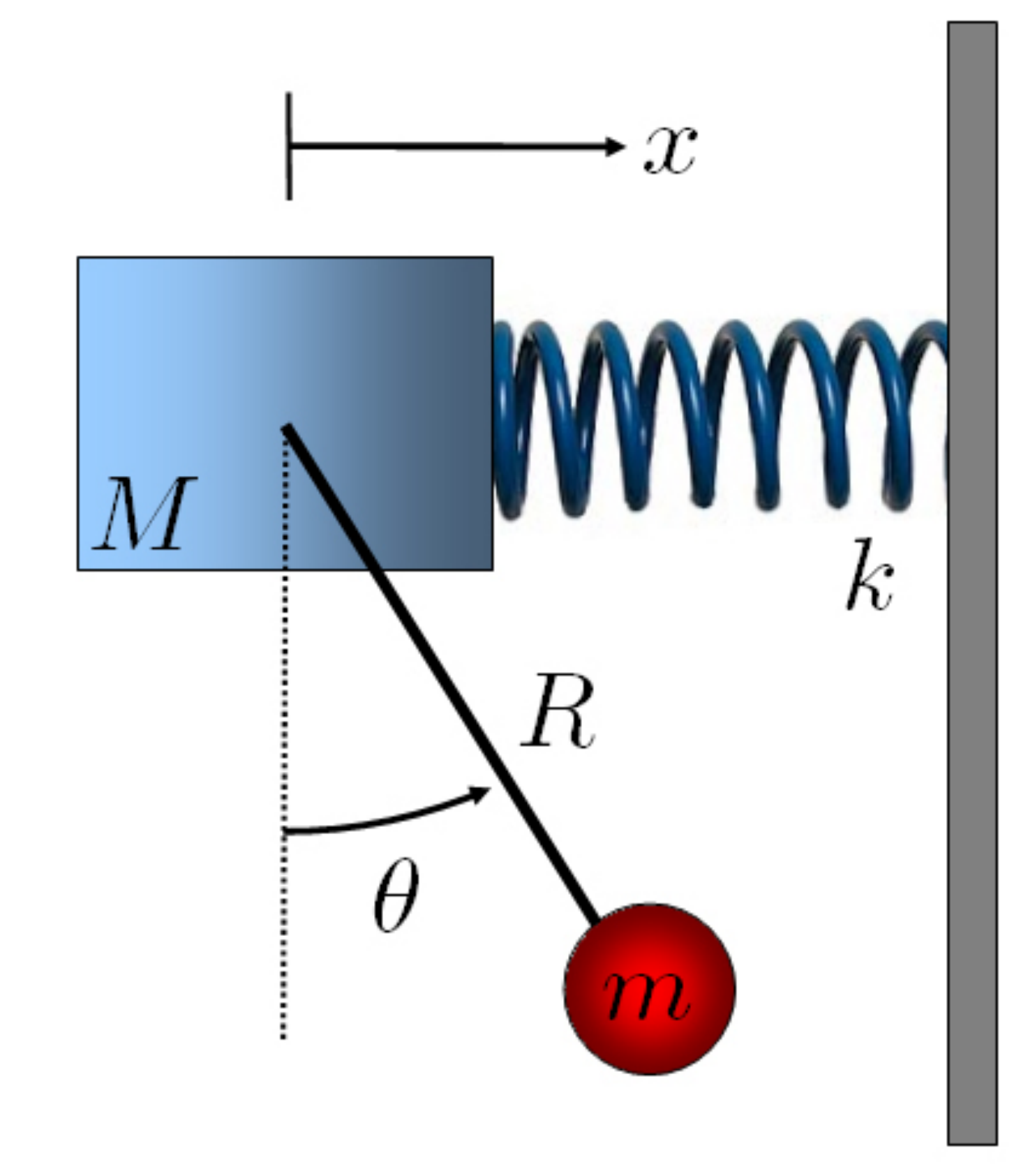}
\caption{A graphical representation of the pendulum on a cart
mounted to a spring.} \label{cartpend}
\end{center}
\end{figure}

\begin{figure*}
\begin{center}
\includegraphics[scale=0.39]{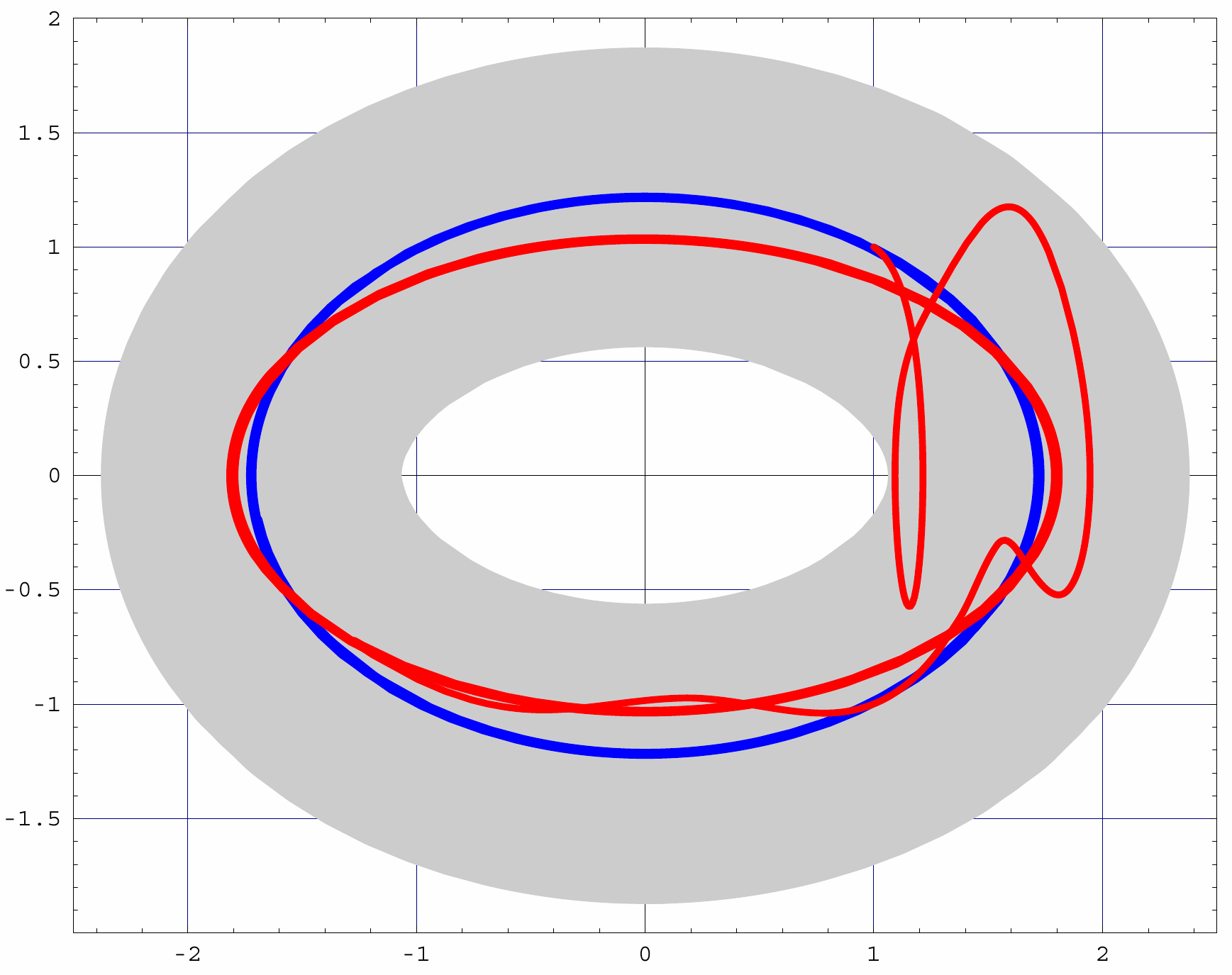}
\includegraphics[scale=0.39]{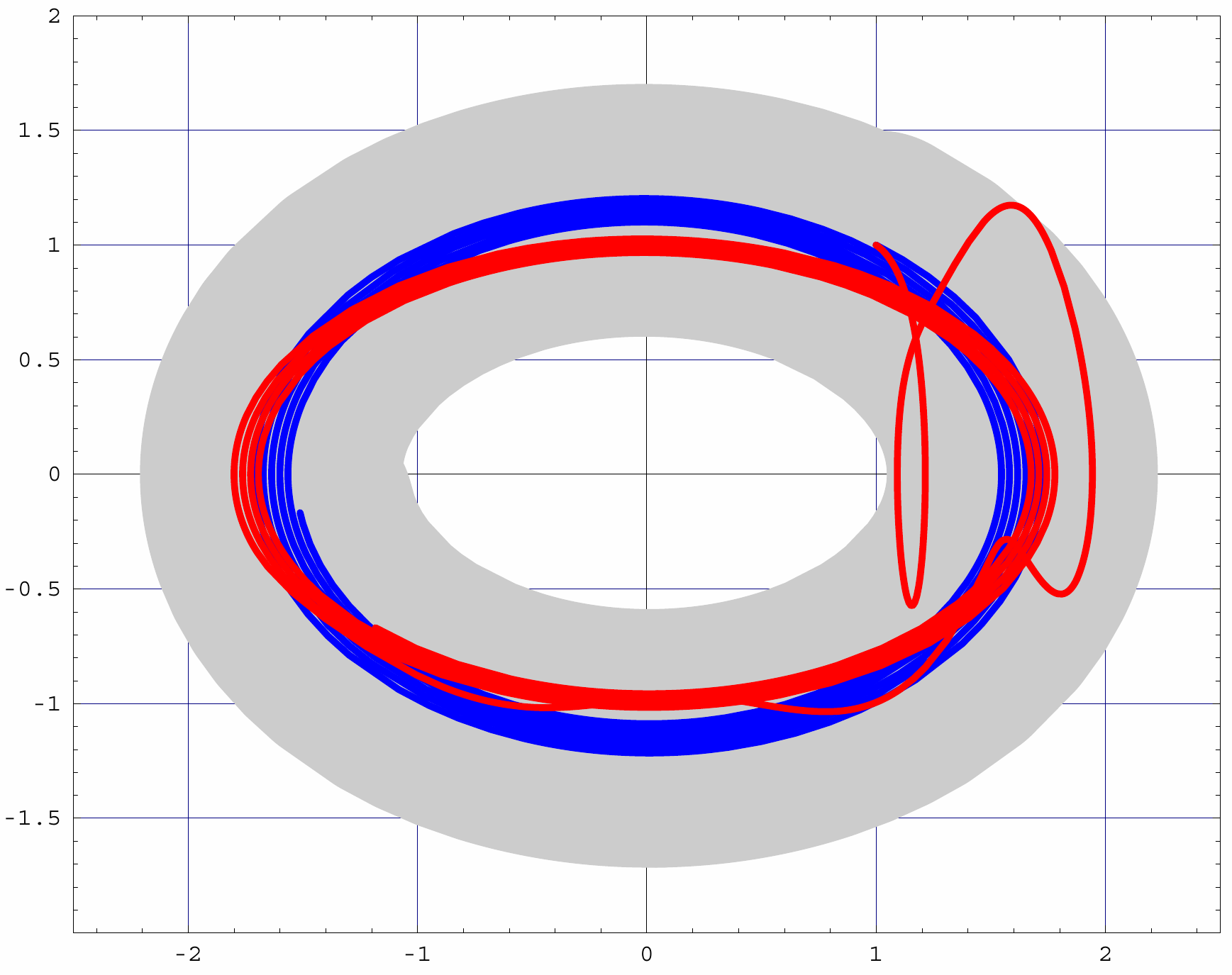}
\includegraphics[scale=0.39]{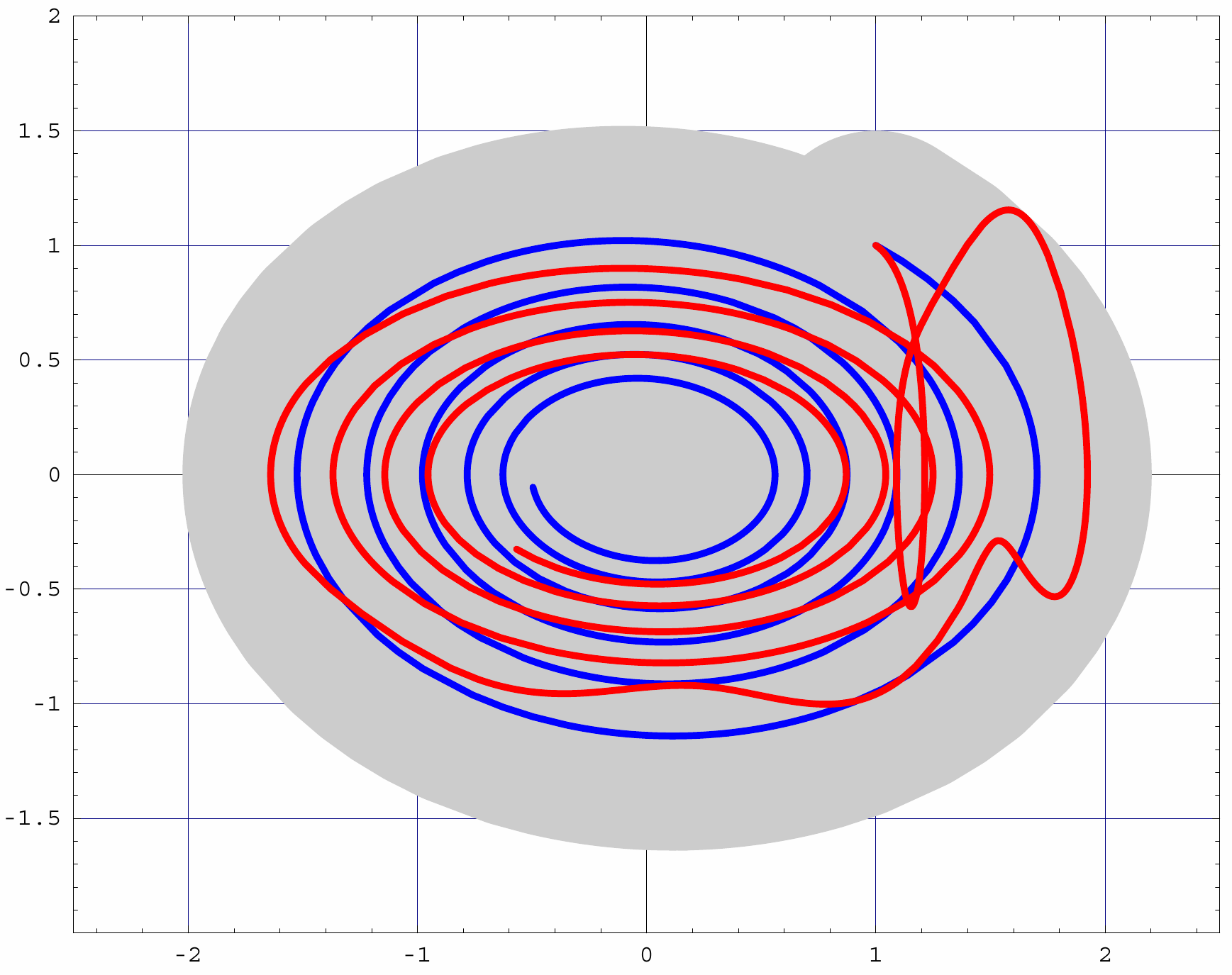}
\includegraphics[scale=0.39]{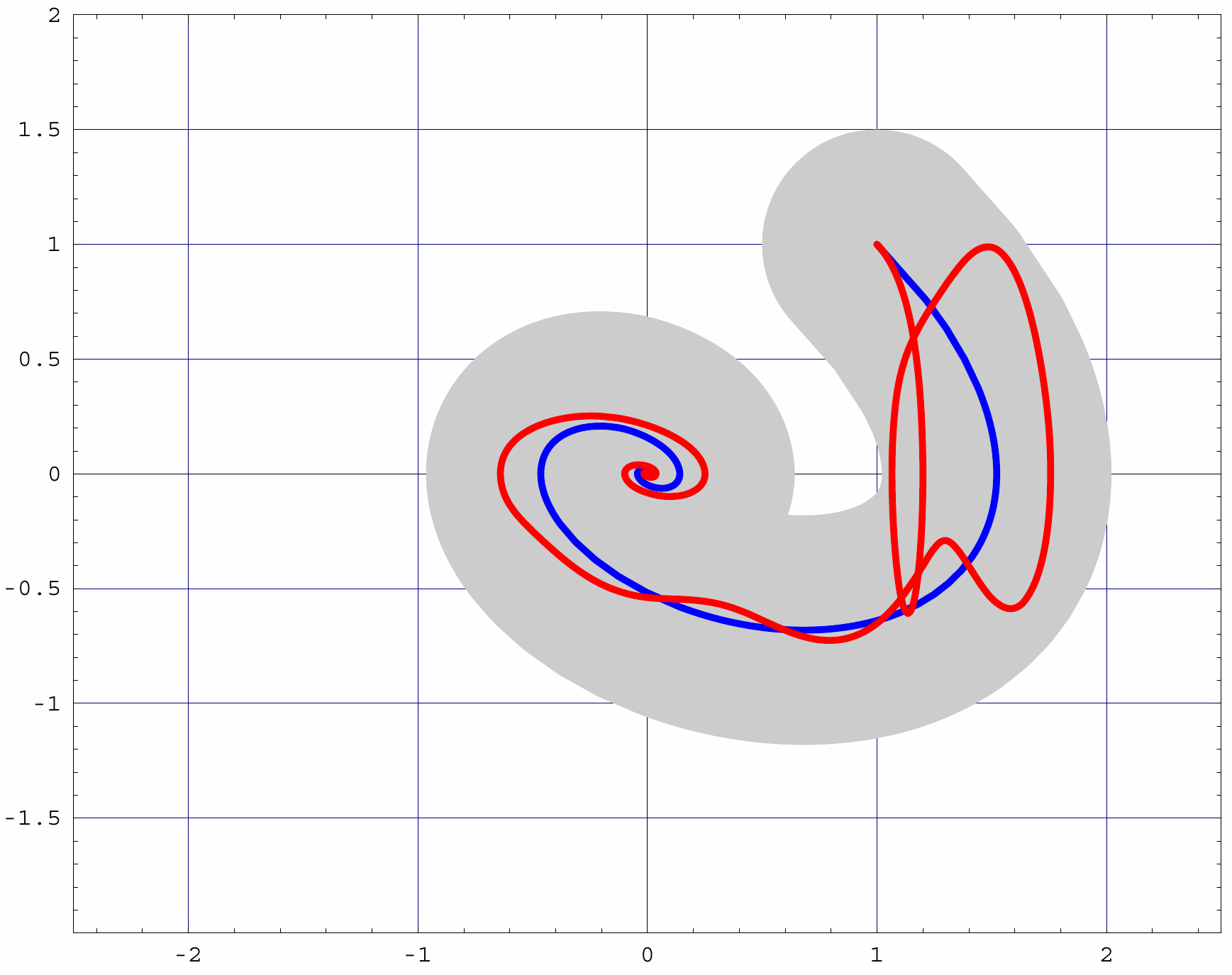}
\caption{A projected trajectory of the full-order system (red) and
a trajectory for the reduced system (blue) for $d =
0.001,0.01,0.1,1$ (from left to right and top to bottom,
respectively).} \label{example2fig}
\end{center}
\end{figure*}

\begin{example}
We now consider a pendulum attached to a cart which is mounted to a
spring (see Figure \ref{cartpend}).  For this example, there are
the following parameters:
\begin{eqnarray}
M & = & \mathrm{mass \:\: of \:\: the \:\: cart}, \nonumber\\
m & = & \mathrm{mass \:\: of \:\: the \:\: pendulum}, \nonumber\\
R & = & \mathrm{length \:\: of \:\: the \:\: rod}, \nonumber\\
k & = & \mathrm{spring \:\: stiffness}, \nonumber\\
g & = & \mathrm{acceleration \:\: due \:\: to \:\: gravity}, \nonumber\\
d & = & \mathrm{friction \:\: constant \:\:for \:\: the \:\:
cart}, \nonumber\\
 b & = & \mathrm{friction \:\: constant
\:\:for \:\: the \:\: pendulum}. \nonumber
\end{eqnarray}

The equations of motion are given by:
\begin{eqnarray}
\label{odeex3}
\dot{x} & = & v\notag\\
\dot{\theta} & = & \omega\notag\\
\dot{v} & = & \frac{1}{M+m\sin^2\theta}\Big(mR\omega^2\sin\theta+m g \sin\theta\cos\theta-kx-dv+\frac{b}{R}\cos\theta\Big)\notag\\
\dot{\omega} & = & \frac{1}{R(M+m\sin^2\theta)}\Big(-mR\omega^2\sin\theta\cos\theta\notag\\
&&-(m+M)g\sin\theta+kx\cos\theta+dv\cos\theta-\Big(1+\frac{M}{m}\Big)\frac{b}{R}\omega\Big)
\end{eqnarray}
where $x$ is the position of the cart, $v$ its velocity, $\theta$
is the angular position of the pendulum and $\omega$ its angular
velocity.

If $\pi_{(x,v)} : \Re^4 \to \Re^2$ is the projection
$\pi_{(x,v)}(x,\theta,v,\omega) = (x,v)$ and $X$ is the vector
field as defined in (\ref{odeex3}), the goal is to reduce $X$ to
$\Re^2$ by eliminating the $\theta$ and $\omega$ variables and thus obtaining $Y$ defined by:
\begin{eqnarray}
\left( \begin{array}{c} \dot{x} \\ \dot{v}
\end{array} \right)  &=& Y(x,v) = T_{(x,\theta,v,\omega)}\pi_{(x,v)} \cdot X(x,0,v,0) =    \left( \begin{array}{c} v \\ - \frac{1}{M} (d v + k x)
\end{array} \right). \nonumber
\end{eqnarray}
The objective is now to show that $X$ and $Y$ are approximately
$\pi_{(x,v)}$-related.  In particular, note that the reduced
system, $Y$, is linear while the full-order system, $X$, is very
nonlinear. This will be discussed in more detail after proving
that they are in fact approximately related.

Stability of $X$, and in particular fiberwise stability, can be proven as in the previous example by noting that $X$ is Hamiltonian for $d=b=0$ and using the Hamiltonian as a Lyapunov function $V$. Consider now the control system:
\begin{eqnarray}
\label{RedC}
F\left((x,v),(\theta,\omega)\right)&=&T\pi_{(x,v)}\cdot X(x,\theta,v,\omega) \\
&=&  \frac{1}{M+m\sin^2\theta}\Big(mR\omega^2\sin\theta -kx\notag\\
&& +m g \sin\theta\cos\theta-dv+\frac{b}{R}\cos\theta\Big)\notag
\end{eqnarray}
with $\theta$ and $\omega$ regarded as inputs. To show that $F$ is
IUBIBSS we first rewrite~(\ref{RedC}) in the form:
\begin{eqnarray}
F\left((x,v),(\theta,\omega)\right)= \frac{1}{M+m}\Big(mR\omega^2\sin\theta
-kx-dv-mR\dot{\omega}\cos\theta\Big) \notag
\end{eqnarray}
and consider the following IISS candidate Lyapunov function:
$$U=\frac{1}{2(m+M)}(x_1-x_2)^2+\frac{1}{2}(v_1-v_2)^2.$$
Its time derivative is given by:
\begin{eqnarray}
\dot{U} & = & -\frac{d}{m+M}(v_1-v_2)^2\notag\\
&&+\frac{mR}{m+M}\Big(\omega_1^2\sin\theta_1-\dot{\omega}_1\cos\theta_1-\omega_2^2\sin\theta_2+\dot{\omega}_2\cos\theta_2\Big)(v_1-v_2).\notag
\end{eqnarray}
Using an argument similar to the one used for the previous
example, we conclude that:
$$\vert v_1-v_2\vert\ge\frac{2mRL}{d}\vert(\theta_1,\omega_1,\dot{\omega}_1)-(\theta_2,\omega_2,\dot{\omega}_2)\vert,$$
with $L$ the Lipschitz constant of the function
$\omega^2\sin\theta-\dot{\omega}\cos\theta$, implies:
$$\dot{U}\le -\frac{d}{2(m+M)}(v_1-v_2)^2,$$
thus showing that $X$ is IISS and in particular also IUBIBSS. That
is, we have established that $X$ and $Y$ are approximately
$\pi_{(x,v)}$-related.

In order to illustrate some of the interesting implications of
approximate reduction, we will compare the reduced system, $Y$,
and the full-order system, $X$, in the case when $R = m = k = b =
1$ and $M = 2$.  It follows that the equations of motion for the
reduced system are given by the linear system:
$$
\left( \begin{array}{c} \dot{x} \\ \dot{v} \end{array} \right) =
 \left( \begin{array}{cc} \phantom{-}0 & \phantom{-}1 \\ -\frac{1}{2}  & \phantom{-}d \\
\end{array} \right) \left( \begin{array}{c} x \\ v \end{array}
\right),
$$
so we can completely characterize the dynamics of the reduced
system: every solution spirals into the origin.  This is in stark
contrast to the dynamics of $X$ (see (\ref{odeex3})) which are very
complex.  The fact that $X$ and $Y$ are approximately related, and
more specifically Theorem \ref{MainTh}, allows us to understand the
dynamics of $X$ through the simple dynamics of $Y$. To be more
specific, because the distance between the projected trajectories of
$X$ and the trajectories of $Y$ is bounded, we know that the
projected trajectories of $X$ will ``essentially'' be spirals.
Moreover, the friction constant $d$ will directly affect the rate of
convergence of these spirals.  Examples of this can be seen in
Figure \ref{example2fig} where $d$ is varied to affect the
convergence of the reduced system, and hence the full order system.
\end{example}

\bibliographystyle{alpha}
\bibliography{/Users/tabuada/Documents/Bibfiles/Control,/Users/tabuada/Documents/Bibfiles/ComputerScience,/Users/tabuada/Documents/Bibfiles/HybridSystems,/Users/tabuada/Documents/Bibfiles/NonLinearControl,/Users/tabuada/Documents/Bibfiles/Algebra,/Users/tabuada/Documents/Bibfiles/DiffGeometry,/Users/tabuada/Documents/Bibfiles/Des,/Users/tabuada/Documents/Bibfiles/GeomMechanics}

\begin{thebibliography}{BKMM96}

\bibitem[AMR88]{MTAA}
R.~Abraham, J.~Marsden, and T.~Ratiu.
\newblock {\em Manifolds, Tensor Analysis and Applications}.
\newblock Applied Mathematical Sciences. Springer-Verlag, 1988.

\bibitem[Ang02]{IncrementalS}
D.~Angeli.
\newblock A {L}yapunov approach to incremental stability properties.
\newblock {\em IEEE Transactions on Automatic Control}, 47(3):410--421, 2002.

\bibitem[ASG00]{SurveyModelReduction}
A.~C. Antoulas, D.~C. Sorensen, and S.~Gugercin.
\newblock A survey of model reduction methods for large-scale systems.
\newblock {\em Contemporary Mathematics}, 280:193--219, 2000.

\bibitem[BDG96]{MReduc}
C.~L. Beck, J.~Doyle, and K.~Glover.
\newblock Model reduction of multidimensional and uncertain systems.
\newblock {\em IEEE Transactions on Automatic Control}, 41(10):1466--1477,
  1996.

\bibitem[BKMM96]{NonH}
A.M. Bloch, P.S. Krishnaprasad, J.E. Marsden, and R.~Murray.
\newblock Nonholonomic mechanical systems with symmetry.
\newblock {\em Archive for Rational Mechanichs and Analysis}, 136(1):21--99,
  1996.

\bibitem[BM00]{UBIBS}
A.~Bacciotti and L.~Mazzi.
\newblock A necessary and sufficient condition for bounded-input bounded-state
  stability of nonlinear systems.
\newblock {\em SIAM Journal on Control and Optimization}, 39(2):478--491, 2000.

\bibitem[GP05]{ApproxBisimDyn}
A.~Girard and G.~J. Pappas.
\newblock Approximate bisimulations for nonlinear dynamical systems.
\newblock In {\em Proceedings of the 44th IEEE Conference on Decision and
  Control}, Seville, Spain, 2005.

\bibitem[GP07]{ApproxBisimLinear}
A.~Girard and G.~J. Pappas.
\newblock Approximate bisimulation relations for constrained linear systems.
\newblock {\em Automatica}, 43(8):1307--1317, 2007.

\bibitem[MR99]{MechSym}
Jerrold~E. Marsden and Tudor~S. Ratiu.
\newblock {\em Introduction to Mechanics and Symmetry}.
\newblock Number~17 in Texts in Applied Mathematics. Springer-Verlag, 2nd
  edition, 1999.

\bibitem[MSVS85]{Dynamics}
Giuseppe Marmo, Alberto Simoni, Bruno Vitale, and Eugene~J. Saletan.
\newblock {\em Dynamical Systems}.
\newblock John Wiley \& Sons, 1985.

\bibitem[MW74]{SymplecticReduction}
J.~E. Marsden and A.~Weinstein.
\newblock Reduction of symplectic manifolds with symmetry.
\newblock {\em Reports on Mathematical Physics}, 5:121--120, 1974.

\bibitem[Son06]{ISS}
E.D. Sontag.
\newblock Input to {S}tate {S}tability: {B}asic concepts and results.
\newblock In P.~Nistri and G.~Stefani, editors, {\em Nonlinear and Optimal
  Control Theory}, pages 163--220. 2006.
\newblock Electronically available at {\tt
  http://www.math.rutgers.edu/$\sim$sontag/}.

\bibitem[Tab07]{Tab07}
P.~Tabuada.
\newblock Approximate simulation relations and finite abstractions of quantized
  control systems.
\newblock In A.~Bemporad, A.~Bicchi, and G.~Buttazzo, editors, {\em Hybrid
  Systems: Computation and Control 2006}, volume 4416 of {\em Lecture Notes in
  Computer Science}, pages 529--542. Springer-Verlag, Pisa, Italy, 2007.

\bibitem[TAJP06]{TAJPCDC07}
P.~Tabuada, A.~Ames, A.~Julius, and G.~J. Pappas.
\newblock Approximate reduction of dynamical systems.
\newblock In {\em Proceedings of the 45th IEEE Conference on Decision and
  Control}, San Diego, CA, December 2006.

\bibitem[TP04]{BisimSCL}
P.~Tabuada and G.~J. Pappas.
\newblock Bisimilar control affine systems.
\newblock {\em Systems and Control Letters}, 52(1):49--58, 2004.

\bibitem[vdS81]{Noether}
A.~van~der Schaft.
\newblock Symmetries and conservation laws for {H}amiltonian systems with
  inputs and outputs: A generalization of {N}oether's theorem.
\newblock {\em Systems and Control Letters}, 1:108--115, 1981.

\bibitem[vdS04]{BisimTAC}
A.~van~der Schaft.
\newblock Equivalence of dynamical systems by bisimulation.
\newblock {\em IEEE Transactions on Automatic Control}, 49(12):2160--2172,
  2004.

\end{thebibliography}

\end{document}